\documentclass[12pt]{article}
\usepackage{amsfonts,amssymb,amsmath,theorem}
\usepackage{oldlfont}
\usepackage{amsfonts}

\oddsidemargin\evensidemargin
\usepackage{array,longtable}
\newcommand{\D}{\displaystyle}
\begin{document}

\begin{center}
    \textbf{  \\[0.5cm]
    MODULAR ESTIMATES IN ORLICZ SPACES \\ AND HAMMERSTEIN OPERATOR EQUATIONS \\ [0.5cm]
    A.V. Lebedev, P.P. Zabreiko} \\[0.5cm]

    \bigskip

    Belarus State University
\end{center}

\vspace{5mm}

\quad  \parbox{14,5cm}{\small \hspace{0.5cm}   Certain
inequalities between the values of the modular and the norm in the
Orlicz spaces are established. These inequalities  are applied
then to the theory of solvability of nonlinear integral equations
of Hammerstein type. \vspace{0.5cm}

\noindent{\bf Keywords:} {\em Orlicz space, modular, Hammerstein
equation, Minty condition, \\ Rothe condition} \vspace{0.5cm}

 }

\bigskip
\bigskip
\bigskip

The aim of the paper is to establish certain inequalities between
the values of the modular and the norm in Orlicz spaces. These
inequalities  turn out to be of use in the theory of nonlinear
integral equations.  \vspace{0.5cm}

{\bf 1.} Let us recall the basic definitions of the theory of
Orlicz spaces  (see, for example,  \cite{KR1}). Henceforth
$(\Omega,{\frak A},\mu)$ is the triple consisting of a set
 $\Omega$,
$\sigma$-algebra  ${\frak A}$ of its subsets and a
$\sigma$-additive measure  $\mu$ defined on  ${\frak A}$. It is
assumed that  $\mu$ is continuous on  $\Omega$  (that is any set
of positive measure can be divided into two sets of equal
measures) and is finite: $\mu(\Omega) < \infty$.

Let  $M : [0, +\infty ) \to [0, +\infty )$ be an $N$-function that
is a convex function  satisfying the conditions
$$
\lim_{u\to 0}\frac{M(u)}{u} = 0, \ \ \ \ \ \lim_{u\to
\infty}\frac{M(u)}{u} = \infty .
$$
On the set of measurable  functions on $\Omega$ one can consider
the functional
$$
{\cal M}(x) := \int\limits_{\Omega} M(|x(\omega)|) \, d\mu (\omega
).
$$
This functional is  conventionally called the  {\it modular}
generated by the $N$-function  $M(\cdot)$. The set  $O_M$ of the
functions $x(\cdot)$ that are measurable on  $\Omega$ and satisfy
the condition ${\cal M}(x) < \infty$ is called the  {\it Orlicz
class}. We shall denote by  $L_M$ the  {\it Orlicz space} that is
the set consisting of the functions $x(\cdot)$ that are measurable
on $\Omega$ and such that for any $x(\cdot)$ there exists
  $\lambda > 0$ such that the following inequality holds
$$
{\cal M}\left(\frac{x}{\lambda}\right) < \infty.
$$
 Orlicz space is a Banach space with respect to each of the
 following two  norms
$$
\|x\|_{\cal L} = \inf \ \left\{\lambda: \ {\cal
M}\left(\frac{x}{\lambda}\right) \le 1\right\}, \qquad \|x\|_{\cal
A} = \inf_{\lambda >0} \ \frac{1 + {\cal M}(\lambda x)}{\lambda}.
$$
 The first of these norms is commonly called the Luxemburg norm
  and the second one --- the
 Orlicz norm
 (in fact both these norms have been  introduced by Orlicz). These norms are equivalent:
\begin{equation}
\label{e1}
\|x\|_{\cal L} \le \|x\|_{\cal A} \le 2\|x\|_{\cal L}.
\end{equation}
The norm calculation in Orlicz spaces even for the most simple
functions is a nontrivial  problem. We present  here the important
formulae for the norms of the characteristic functions
$\chi_D(\cdot) \ (D \in {\frak A})$:
$$
\|\chi_D\|_{\cal L} = \frac1{M^{-1}\left(\D\frac1{\mu(D)}\right)},
\qquad \|\chi_D\|_{\cal A} = \mu(D) \,
(M^*)^{-1}\left(\frac1{\mu(D)}\right),
$$
 where  $M^*(\cdot)$ is the  $N$-function  {\it dual to the  $N$-function
$M(\cdot)$}, that is the function given by the equality
 $M^*(u) := \sup \ \{uv - M(v): 0 \leq v <
\infty\}$.

An important role in the theory of Orlicz spaces is played by the
subspace  $L_M^\circ$ of measurable functions $x(\cdot)$
 on  $\Omega$  such that for each
 $\lambda
> 0$ the following inequality takes place
$$
{\cal M}\left(\frac{x}{\lambda}\right) < \infty.
$$
This subspace coincides with the closure of the set of bounded
functions in  Orlicz space. The following embeddings are true
$L_M^\circ \subseteq O_M \subseteq L_M$. Under the presupposed
assumptions on  $(\Omega,{\frak A},\mu)$ each of the equalities
$L_M^\circ = O_M$ and $O_M = L_M$ is equivalent to the  {\it
$\Delta_2$-condition}:
$$
\mathop{\overline{\rm lim}}\limits_{u \to \infty} \
\frac{M(2u)}{M(u)} < \infty.
$$
 In the general case the following embeddings take place
$$
\{x \in L_M: \ d(x,L_M^\circ) < 1\} \subset O_M \subset \{x \in
L_M: \ d(x,L_M^\circ) \le 1\},
$$
 where  $d(\cdot,L_M^o)$ is the distance from the corresponding element to
  $L_M^\circ$ (this distance is the same in the both norms!).

Let us also note the fundamental equalities:
$$
((L_M^\circ)_{\cal L})^* = ({L_{M^*}})_{\cal A}, \qquad
((L_M^\circ)_{\cal A})^* = ({L_{M^*}})_{\cal L},
$$
where by the symbols  $_{\cal L}$ and  $_{\cal A}$ we mark  the
spaces equipped with the corresponding norms. These equalities
mean in particular that under the passage to the dual spaces the
 Luxemburg and Orlicz norms interchange their places.

In what follows the notation  $\langle \cdot,\cdot \rangle$ will
mean the standard coupling of the spaces    $L_{M^*}$ and $L_M$:
$$
\langle y,x \rangle = \int\limits_\Omega
y(\omega)\overline{x(\omega)} \, d\mu(\omega).
$$
One can verify  that
$$
\langle y,x \rangle \le  \|y\|_{({L_{M^*}})_{\cal L}}
\|x\|_{(L_M)_{\cal A}}
$$
 and
 $$
\langle y,x \rangle \le  \|y\|_{({L_{M^*}})_{\cal A}}
\|x\|_{(L_M)_{\cal L}}.
 $$
This along with  (\ref{e1}) implies
$$
\langle y,x \rangle \le  \|y\|_{({L_{M^*}})_{\cal A}}
\|x\|_{(L_M)_{\cal A}}
$$
and
$$
\langle y,x \rangle \le  2\|y\|_{({L_{M^*}})_{\cal L}}
\|x\|_{(L_M)_{\cal L}}.
$$

As the particular examples of Orlicz spaces we have the classical
Lebesgue spaces $L_p \ (1 < p < \infty)$ that correspond to the
$N$-functions   $M(u) = |u|^p \ (1 < p < \infty)$.

The modular  is linked  with  the Luxemburg norm  by the relations
\begin{equation}
\label{(1)}
\|x\|_{\cal L} \le 1 \ \ \Rightarrow \ \ {\cal M}(x) \le
\|x\|_{\cal L}, \qquad \|x\|_{\cal L} \ge 1 \ \ \Rightarrow \ \
{\cal M}(x) \ge \|x\|_{\cal L};
\end{equation}
and with the  Orlicz norm by the inequality
\begin{equation}
\label{(2)}
{\cal M}(\lambda x) \ge \lambda\|x\|_{\cal A} - 1 , \qquad 0 <
\lambda < \infty.
\end{equation}

It was observed in  \cite{G}  that the next  equality
\begin{equation}
\label{(3)}
\lim_{\|x\| \to \infty} \ \frac{{\cal M}(x)}{\|x\|} =
\infty
\end{equation}
plays a significant role in applications.
 As it was shown in  \cite{R} this equality  is not always true.
 The necessary and sufficient condition for its validity
 is the equality
 \begin{equation}
\label{(4)}
 \lim_{u \to \infty} \ \frac{M(ku)}{u(M^*)^{-1}(u)} = \infty ,\qquad
1 < k < \infty .
\end{equation}
  As the examples
 of  $N$-functions that do not satisfy condition  (\ref{(4)})
 one can take the functions  $(1 + u) \ln (1 + u) -
u$ \, and \, $u \sqrt{\ln (1 + u)}$
(see \cite{R}). These functions are characterized
by the property that the corresponding dual functions satisfy the
so called
 $\Delta_3$-condition  (that is  for some  $k > 1$ and large  $u$
 the inequality  $uM(u)
\le M(ku)$ holds; see \cite{KR1}). We also remark  that a number
of properties of Orlicz space associated with relation  (\ref{(3)}) were
considered in  \cite{H}.

However in the particular case  $M(u) = u^p \ \ (1 < p < \infty)$
the Luxemburg norm  coincides with the standard norm on $L_p$ and
satisfies  the equality
\begin{equation}
\label{(4*)}
{\cal M}(x) = \|x\|^p,
\end{equation}
 This equality is stronger than   relations
 (\ref{(1)}) and   (\ref{(3)}).
Thus one naturally arrives at the problem of a possibility of
refinement of these relations   for arbitrary Orlicz spaces.

Henceforth we confine ourselves to the estimates in the Luxemburg
norm. Therefore for the sake of brevity  the sign $_{\cal
L}$ in the notation of the norm $\|\cdot\|_{\cal L}$ is omitted.
\vspace{0.5cm}

{\bf 2.} In this section we  discuss the problem of
validity  in Orlicz spaces of the estimates of the type
\begin{equation}
\label{(5)}
{\cal M}(x) \ge \phi(\|x\|), \qquad \|x\| \ge R ,
\end{equation}
where  $\phi(\cdot)$ is a certain positive function and  $R$ is a sufficiently
large number. Note that  (\ref{(1)}) implies that one can always take as the function  $\phi(\cdot)$   the function
\begin{equation}
\label{(6)}
\phi(\lambda) = \left\{\begin{array}{ccc}0 & \text{при} & 0 \le \lambda < 1 \\
\lambda & \text{при} & 1 \le \lambda < \infty \end{array}
\right.
\end{equation}
 But since we are aimed at (\ref{(3)})  the functions $\phi(\cdot)$ of prime interest are those
 who grow faster than linear functions at infinity.

The next statement gives a partial answer to the question.

{\bf Theorem  1}. {\it Let  an $N$-function  $M(\cdot)$ satisfy
the inequality
\begin{equation}
\label{(7)}
M(\lambda u) \ge \phi(\lambda) \cdot M(u), \qquad  \lambda \geq 1
,
\end{equation}
 where  $\phi(\cdot)$ is a  positive function.
Then
\begin{equation}
\label{(8)}
{\cal M}(x) \ge \phi(\|x\|), \qquad \|x\| \ge 1 .
\end{equation}}

$\square$ The proof is quite simple.
 Indeed, let $\|x\|
\ge 1$.  Then we have by (\ref{(7)}) and (\ref{(1)})
$$
{\cal M}(x) = {\cal M}\left({\|x\|} \, \frac{x}{\|x\|}\right) \ge
\phi\left({\|x\|}\right) \ {\cal M}\left(\frac{x}{\|x\|}\right)
\ge \phi\left({\|x\|}\right). \ \ \ \ \blacksquare
$$

As the examples of application of this theorem one can consider not only the
 $N$-functions  $M(u) = u^p \ (1 < p < \infty)$
corresponding to the classical Lebesgue spaces  $L_p$ (for these functions one can take
$\phi(\lambda) = \lambda^p$), but also the  $N$-functions  $e^{u}
- u - 1$ and  $e^{u^2} - 1$; for these functions inequality  (\ref{(7)})
is satisfied with $\phi(\lambda) = \lambda^2$. Thus Theorem 1
implies that  in the corresponding spaces the  Luxemburg norm
satisfies the inequality
$$
{\cal M}(x) \ge \|x\|^2 , \qquad \|x\| \ge 1.
$$

For the functions   $e^{u} - u - 1$ and  $e^{u^2} - 1$ mentioned
above it is natural to expect that the inequality of the form  (\ref{(5)})
for large  $\|x\|$ is satisfied with a function $\phi(\cdot)$ that
grows at infinity essentially faster than $\lambda^2$. It turns
out that this is really true but Theorem 1 is not enough to prove
this fact: by means of this theorem the function  $\lambda^2$ for
both  the functions $e^{u} - u - 1$ and  $e^{u^2} - 1$ cannot be
changed for the one that grows faster.

{\bf Theorem 2}. {\it Let an  $N$-function  $M(\cdot)$ satisfy the
inequality
\begin{equation}
\label{(9)}
M(\lambda u) \ge \phi(\lambda , u) \cdot M(u), \qquad  \lambda \ge
1,
\end{equation}
where  $\phi(\cdot,\cdot)$ is a  positive and
nondecreasing  with respect to the first and the second argument
function.
 Then    the following inequality  holds
\begin{equation}
\label{(10)}
{\cal M}(x) \ge \phi\left(\frac{\| x\|}{R}, \frac{R-1}{R \| {\bf
1}\|}\right) \qquad  {{for}}\ \ \|x\| \ge R > 1 ,
\end{equation}
where ${\bf 1}$ is the function identically equal to $1$.}

$\square$ Let a function  $x \in L_M$ satisfy the condition $\|x\|
\ge 1$ and in addition for each point  $\omega \in \Omega$ where
$x(\omega) \ne 0$ we have  $|x(\omega)| \ge h\|x\|$, where $h$ is
a certain positive number. Repeating word by word the argument of
the proof of Theorem 1 and using (\ref{(9)}) instead of  (\ref{(7)}) and the
monotone property of $\phi$ with respect to $u$ we get
\begin{equation}
\label{(11)}
{\cal M}(x) \ge \phi(\|x\|, h).
\end{equation}

Now let  $x \in L_M$ be a function for which  $\|x\| \ge R>1$. Set
$$
h := \frac{R-1}{R \| {\bf 1}\|}
$$
 and
$$
x_h(\omega) := \left\{\begin{array}{ccc}
x(\omega), & \text{при} & |x(\omega)| \ge h\|x\|, \\
0, & \text{при} & |x(\omega)| < h\|x\|.
\end{array}\right.
$$
 Then
$$
\|x_h\| \ge \|x\| - \|x - x_h\| \ge \|x\| - h\|{\bf 1}\| \cdot
\|x\| = \frac{\|x\|}{R}\ge 1.
$$
 Now  from the inequality  $|x(\omega)| \ge |x_h(\omega)|$,
 (\ref{(11)}) and the monotone property of $\phi$ with respect to $\lambda$
we get
$$
{\cal M}(x) \ge {\cal M}(x_h) \ge \phi\left(\frac{\|x\|}{R},
h\right) \ \ {\rm for}\ \ \|x\| \ge R>1 . \ \ \blacksquare
$$
{\bf Remark.} Clearly if we additionally presuppose that $\phi$ is continuous then in
(\ref{(10)}) we can also take $R=1$.

Let us consider as an example  of  $N$-function  the function $e^{u} - u - 1$.
Evidently the function
\begin{equation}
\label{(12)}
\phi(\lambda,u) = \frac{e^{\lambda u} - \lambda u - 1}{e^u - u -
1}
\end{equation}
 is nondecreasing with respect to $\lambda $. Observe that it is also
 nondecreasing with respect to $u$ on  $[0,\infty)$ for  $\lambda \ge 1$.
 Indeed, it is enough to
 verify the nonnegativity of its derivative; but the latter is equivalent to the
 inequality
$$
(\lambda e^{\lambda u} - \lambda)(e^u - u - 1) - (e^{\lambda u} -
\lambda u - 1)(e^u - 1) \ge 0
$$
which  is equivalent to
$$
\left(\sum_{j=2}^\infty \frac{j\lambda^j u^{j-1}}{j!}\right)
\left(\sum_{k=2}^\infty \frac{u^k}{k!}\right) -
\left(\sum_{j=2}^\infty \frac{\lambda^j u^j}{j!}\right)
\left(\sum_{k=2}^\infty \frac{ku^{k-1}}{k!}\right) \ge 0.
$$
Collecting similar terms in the  left hand part of this inequality
we obtain
$$
  \sum_{2 \le k < j < \infty}
\frac{(j - k)(\lambda^j - \lambda^k)u^{j+k-1}}{j!k!} \ge 0.
$$
The latter inequality is evident since  $\lambda \ge 1$.

In addition we have
$$
e^{\lambda u} - \lambda u - 1 = \phi(\lambda, u) \ (e^{u} - u - 1).
$$
 Now Theorem  2 implies the validity  of estimate
  (\ref{(10)}) with  function  (\ref{(12)}).

In an analogous way one can consider the  $N$-function  $e^{u^2} -
1$ taking here
$$
\phi(\lambda, u) := \frac{e^{\lambda^2u^2} - 1}{e^{u^2} - 1} \ \ .
$$

It can be shown that Theorems  1 and 2 contain the results of Ja. B.
Rutitzky from \cite{R}. Note that the statement of the main
theorem of the latter paper contains an  inexactness: condition
(\ref{(5)}) in the proof of sufficiency is used not for large $u$ but for
all $u$. \vspace{0.5cm}

{\bf 3.} In the paper  \cite{ZR} there was investigated in Orlicz
spaces the nonlinear integral Hammerstein equation of the form
\begin{equation}
\label{(13)}
x = Sf(x) + g,
\end{equation}
where  $f$ is a nonlinear superposition operator
$$
fx(\omega) = f(\omega,x(\omega)),
$$
where $f(\omega,u): \ \Omega \times {\Bbb R} \to {\Bbb R}$ is a
function satisfying the Caratheodory conditions, $S$ is a linear
operator and  $g$ is a known function. As  examples  of (\ref{(13)})
 one can consider the
nonlinear singular integral equations  (see, for example,
 \cite{GM,ZR}); in these examples  $\Omega \subseteq {\Bbb
R}^2$ is a certain  (open or closed) sufficiently smooth curve
 and  $S$ is the linear integral operator with the Cauchy type kernel.

Unfortunately in   \cite{ZR} the signs in a number of inequalities
happened to be mixed up and as a result the statements of Lemma 4
and Theorem  2 in \cite{ZR} on the conditions of existence and
uniqueness of the solution to equation
 (\ref{(13)}) are false. Hereafter we give the general scheme
 of the investigation suggested in  \cite{ZR} with the necessary refinement.
 A number of new important additional
 observations are presented as well.

To start with we describe the general scheme of the proof of {\it
existence theorems}. Henceforth in this Section we shall consider
the situation $L_M \subset L_2 \subset L_{M^*}$ that implies in
particular the inequality
\begin{equation}
\label{l2}
\gamma \|x\| \ge\|x\|_2,
\end{equation}
where $\|x\|_2$ is the $L_2$-norm of $x$ and $\gamma >0$ is a
certain constant. Let the superposition operator $f$ act from an
Orlicz space $L_M$ into its dual Orlicz space $L_{M^*}$ (recall
that $M^*(\cdot)$ is the $N$-function dual to the $N$-function
$M(\cdot)$). Further, let the operator $S$ be defined on $L_{M^*}$
and take values in the space of measurable functions, in addition
let there exist a linear operator $T$ acting from  $L_M$ into
$L_{M^*}$ such that
\begin{equation}
\label{(14)}
STx = x, \qquad x \in L_M .
\end{equation}
Finally, let the function  $g$ also belong to   $L_M$.
 Under these assumptions to prove the solvability of equation
(\ref{(13)}) one can consider an auxiliary equation   $\Phi x = 0$, where
$$
\Phi x = Tx - f(x) - Tg.
$$
 Indeed, applying the operator $S$ to this equation
 we arrive at equation (\ref{(13)}) which means that every
 solution  $x \in L_M$
to the equation  $\Phi x = 0$ is also the solution in  $L_M$ to
 equation  (\ref{(13)}).

Under a number of  natural constraints on the operators  $S$ and
$f$ the operator  $\Phi$ as the operator from  $L_M$ into
$L_{M^*}$ turns out to be monotone in the sense of Minty:
$$
{\rm Re} \, \langle \Phi x_1 - \Phi x_2,x_1 - x_2 \rangle \ \ge 0
, \qquad x_1, x_2 \in L_M,
$$
 and on the balls  $\{x: \ \|x\| \le R\}$ of the space  $L_M$ it satisfies the
 Rothe condition
$$
{\rm Re} \, \langle \Phi x,x \rangle \ \ge 0 , \qquad \|x\| = R.
$$

We recall that  the Minty monotonicity  along with the Rothe
condition imply the existence of a solution to the equation $\Phi
x = 0$ (see \cite{ZR}) and therefore
 the existence of a solution to (\ref{(13)}) as well.

Now we shall describe  the corresponding constrains. Let the
linear operator $T$ satisfies the inequality
\begin{equation}
\label{(15)}
{\rm Re} \, \langle Tx,x \rangle \ge \sigma \, \|x\|_2^2 ,\qquad x
\in L_M .
\end{equation}
 Then, if the function
$f(\omega,u)$ satisfies with respect to the
 argument  $u$ the inequality
\begin{equation}
\label{(16)}
{\rm Re} \, (f(\omega,u_1) - f(\omega,u_2))(\overline{u}_1 -
\overline{u}_2) \le \delta |u_1 - u_2|^2 , \qquad \omega \in
\Omega, \ u_1, u_2 \in {\Bbb C}
\end{equation}
 (this inequality means that the function  $\delta u - f(\omega,u)$
is monotone with respect to the argument  $u$) then the operator
 $\Phi$ is evidently monotone in the sense of Minty on  $L_M$ provided
$\sigma \ge \delta $.

Further, if the function  $f(\omega,u)$ satisfies the inequality
\begin{equation}
\label{eq}
- \, {\rm Re} \, \overline{u}f(\omega,u) \ge aM(b|u|) - c(\omega),
\qquad \omega \in \Omega, \ u \in {\Bbb C},
\end{equation}
 where  $a, b > 0$ and  $c(\cdot)$ is an integrable on  $\Omega$ function then
 for a certain function  $\phi(\cdot)$ (see (\ref{(1)}), (\ref{(6)}) and Theorems  1 and 2) the operator
$f$ satisfies the inequality of the form
$$
- \, {\rm Re} \, \langle fx,x \rangle \ \ge a \phi(b\|x\|) - c ,
\qquad \|x\| \ge R;
$$
here  $c = \|c(\cdot)\|_1$ is the $L_1$-norm of $c( \cdot)$.
Therefore the operator $\Phi$ satisfies the estimate
\begin{equation}
\label{(17)}
{\rm Re}\,\langle \Phi x,x \rangle \ge \sigma \|x\|_2^2 + a
\phi(b\|x\|) - 2\|Tg\|_{L_{M^*}} \, \|x\| - c , \qquad \|x\| \ge
R .
\end{equation}
Observe  that
$$
 \|x\|^{-1}\left(\sigma \|x\|_2^2 + a
\phi(b\|x\|) - 2\|Tg\|_{L_{M^*}} \, \|x\| - c \right)=$$
\begin{equation}
\label{*17} \frac{\sigma \|x\|_2^2}{\|x\|} +
ab\frac{\phi(b\|x\|)}{b\|x\|} - 2\|Tg\|_{L_{M^*}} - \|x\|^{-1}\, c
\end{equation}

Therefore we get  the following statement:
{\em let  $\sigma \ge 0$ and  $r^{-1}\phi(r) \to \infty$ as $r \to
\infty$ then the Rothe condition is satisfied for $R$ large
enough.}

 In other words in
this case it should  be fulfilled  condition (\ref{(3)})
considered in \cite{G} and \cite{R}, that is in this case an
arbitrary Orlicz space can not be used  (however if the norm
$\|Tg\|_{L_{M^*}}$    is sufficiently small one can take as the
function  $\phi(\cdot)$  in (\ref{*17}) function  (\ref{(6)}), see
the left hand inequality in  (\ref{e1}); this means that equation
(\ref{(13)}) with  $g$ of this type can be considered in any of
the Orlicz spaces).

Now let $\sigma < 0$. Recalling that  $\gamma \|x\| \ge\|x\|_2$
(see (\ref{l2}))
 we conclude that
 equality (\ref{*17}) implies the inequality
$$
 \|x\|^{-1}\left(\sigma \|x\|_2^2 + a
\phi(b\|x\|) - 2\|Tg\|_{L_{M^*}} \, \|x\| - c \right)\geq$$
\begin{equation}
\label{*18} {\sigma \gamma^2 \|x\|} +
ab\frac{\phi(b\|x\|)}{b\|x\|} - 2\|Tg\|_{L_{M^*}} - \|x\|^{-1}\, c
\end{equation}

Thus in this situation  we get the following statement: {\em let $\sigma <0$ then the Rothe condition is satisfied for $R$
large enough if $\phi (\cdot)$ has greater than  quadric rate of
growth at infinity.}

Therefore in this case  an arbitrary Orlicz space is not suitable.
Theorems 1 and 2 give us a possibility to indicate the conditions
on the $N$-function $M(\cdot)$ under which the function
$\phi(\cdot)$ of this sort does exist.

Now let us describe the general scheme of the proof of   {\it
uniqueness theorems} for equation  (\ref{(13)}). Suppose
that for the operator  $S$ there exists an operator  $T$ having
the property
$$
TSz = z ,\qquad z \in L_{M^*}.
$$
The equalities  $x_1 = Sfx_1 +g$ and  $x_2 = Sfx_2 +g$ for  $x_1,
x_2 \in L_M$ imply  $x_1 - x_2 = S(fx_1 - fx_2)$. Applying the
operator $T$  to this equality we get $T(x_1 - x_2) = (fx_1 -
fx_2)$. From this under the fulfilment of (\ref{(15)}) and
(\ref{(16)}) we obtain
$$
0={\rm Re} \, \langle (T(x_1 - x_2) - (fx_1 - fx_2)),(x_1 - x_2)
\rangle \ \ge (\sigma - \delta) \|x_1 - x_2\|_2^2,
$$
 and thus if   $\sigma - \delta > 0$ we have  $x_1 = x_2$.
\vspace{0.5cm}

{\bf 4.} Estimates (\ref{(5)}) deduced in Section  {\bf 2} were
true only for sufficiently large   $\| x\|$. For small  $\| x\|$
the similar inequalities are not valid in the general case. In
particular one can easily show that for  $N$-functions  $M(\cdot)$
satisfying the relation
\begin{equation}
\label{(18)}
\mathop{\underline{\rm lim}}\limits_{u \to \infty} \
\frac{M(\lambda u)}{M(u)} = 0 ,\qquad 0 < \lambda < 1,
\end{equation}
the next equality is true:
$$
\inf \ \{{\cal M}(x): \ \ \|x\| \ge R\} = 0, \ \  \ \ 0 < R < 1
$$
 This equality means that the inequalities of the form  (\ref{(5)}) with small  $\| x\|$
and positive functions  $\phi(\cdot)$ are impossible for
$N$-functions satisfying (\ref{(18)}). This remark is true in
particular for the functions $e^{u} - u - 1$ and $e^{u^2} - 1$
considered above.

The foregoing observation implies that the estimate we are
interested in is possible only if
$$
\mathop{\underline{\rm lim}}\limits_{u \to \infty} \
\frac{M(\lambda u)}{M(u)} > 0,  \qquad 0 < \lambda < 1.
$$
 It is easy to see that the latter condition
 is equivalent to the condition that
 the
$N$-function  $M(u)$ satisfies  $\Delta_2$-condition. It turns out
that  $\Delta_2$-condition guarantees the satisfaction of
estimates  (\ref{(5)}) for   $\| x\|\le 1$ and $\phi$ taken from
(\ref{(6)}).

 {\bf Theorem  3}. {\it Let an  $N$-function
$M(\cdot)$ satisfy the inequality
\begin{equation}
\label{(19)}
M(\lambda u) \ge \phi(\lambda) \cdot M(u), \qquad 0 < \lambda \le
1,
\end{equation}
 where  $\phi(\cdot)$ is a  positive function.
Then the next inequality is true
\begin{equation}
\label{(20)}
{\cal M}(x) \ge \phi(\|x\|), \qquad \|x\| \le 1.
\end{equation}}

$\square$ The proof is the word by word repetition of the proof of
Theorem 1. $\blacksquare$

As the examples of  $N$-functions satisfying inequality
(\ref{(19)}) one can consider the functions  $M_1 =(1 + u)\ln (1 +
u) - u$ and $M_2=u^p \ln (1 + u) \ (1 < p < \infty)$. For the
first of these functions the corresponding function
$\phi(\cdot)$ is defined by the equality
 $\phi(\lambda) = \lambda^2$. Indeed, for each  $u > 0$ and  $0 < \lambda < 1$
 by the Cauchy  theorem we have
$$
\frac{M_1(\lambda u)}{M_1(u)} = \frac{\lambda M_1'(\lambda
u_1)}{M_1'(u_1)} = \frac{\lambda^2 M_1''(\lambda
u_2)}{M_1''(u_2)}, \qquad 0 < u_2 < u_1 < u.
$$
 Since  $M_1''(u) = \D\frac1{1 + u}$ is monotone decreasing function on  $(0,\infty)$
it follows from the latter equality that  $M_1(\lambda u) \ge
\lambda^2 M_1(u)$. This inequality cannot be refined since by the
L'Hospital theorem we have
$$
\lim_{u \to 0} \ \frac{M_1(\lambda u)}{M_1(u)} = \lim_{u \to 0} \
\frac{\lambda\ln (1 + \lambda u)}{\ln (1 + u)} = \lim_{u \to 0} \
\frac{\lambda^2(1 + u)}{1 + \lambda u} = \lambda^2.
$$
 Observe also that the equalities
$$
\lim_{u \to \infty} \ \frac{M_1(\lambda u)}{M_1(u)} = \lim_{u \to
\infty} \ \frac{\lambda\ln (1 + \lambda u)}{\ln (1 + u)} =
\frac{\lambda^2(1 + u)}{1 + \lambda u} = \lambda
$$
and the monotone property of the function
$$
\frac{\lambda^2(1
+ u)}{1 + \lambda u}, \ \ u \in (0, \infty )
$$
show that the ratio between  $M_1(\lambda u)$ and $M_1(u)$ for  $0
< \lambda < 1$ is contained in the interval $(\lambda^2,\lambda)$.
Thus $M_1$ satisfies the following inequalities
\begin{equation}
\label{M1} \lambda^2 M_1 (u)\le M_1 (\lambda u)\le \lambda M_1(u),
\ \ \ \ 0 < \lambda < 1.
\end{equation}
\smallskip

For the  function $M_2$ we take   $\phi_2(\lambda) =
\lambda^{p+1}$. Since for  $u > 0$ we have
\begin{equation}
\label{eM2}
\frac{(\lambda u)^p \ln (1 + \lambda u)}{u^p \ln (1 +
u)} = \lambda^p \ \frac{\ln (1 + \lambda u)}{\ln (1 + u)} =
\lambda^{p+1} \ \frac{1 + u_1}{1 + \lambda u_1} \ge \lambda^{p+1},
\end{equation}
where by the Cauchy theorem   $0 < u_1 < u$. The equalities
$$\lim_{u \to 0} \ \frac{M_2(\lambda u)}{M_2(u)} = \lim_{u \to 0} \
\lambda^p \ \frac{\ln (1 + \lambda u)}{\ln (1 + u)} =
\lambda^{p+1} \ \lim_{u \to 0} \ \frac{1 + u}{1 + \lambda u} =
\lambda^{p+1}
$$
 show that the inequality obtained can not be refined.

 Observe also that relations (\ref{eM2}) and the inequality
 $$
\frac{(\lambda u)^p \ln (1+\lambda u)}{u^p \ln (1+u)} \le
\lambda^p, \ \ \ \ \  0< \lambda <1 ,
 $$
show that for this function the ratio between   $M(\lambda u)$ and
$M(u)$ for  $0 < \lambda < 1$ is contained in the interval
$(\lambda^{p+1},\lambda^p)$.

Thus $M_2$ satisfies the following inequalities
\begin{equation}
\label{M2} \lambda^{p+1}M_2 (u)\le M_2 (\lambda u)\le \lambda^p
M_2 (u), \ \ \ \ 0 < \lambda < 1.
\end{equation}

Inequalities (\ref{M2}) and (\ref{M1}) show in particular that if
we take $M_3 =M_1 + M_2$ then we get the $N$-function satisfying
the inequalities
\begin{equation}
\label{M3} \lambda^{p+1}M_3 (u)\le M_3 (\lambda u)\le \lambda M_3
(u), \ \ \ \ 0 < \lambda < 1.
\end{equation}

Just as in the case of   Theorem  1 inequality  (\ref{(19)}) (the
analogue to inequality
 (\ref{(7)})) is too restrictive. Let us present the analogue to Theorem 2.

{\bf Theorem  4}. {\it Let an  $N$-function  $M(\cdot)$ satisfy
the inequality
\begin{equation}
\label{(21)}
M(\lambda u) \ge \phi(\lambda, u) \cdot M(u), \qquad 0 < \lambda <
1,
\end{equation}
 where $\phi(\cdot,\cdot)$ is a  positive and
 nondecreasing with respect to the first and the second argument function.
 Then the following inequalities hold
\begin{equation}
\label{(22)} {\cal M}(x) \ge \phi((1 - h\|{\bf 1}\|)\|x\|, h),
\qquad 0 \ne\|x\| \leq 1, \ 0 < h < \frac1{\|{\bf 1}\|}.
 \end{equation}}

$\square$ Let a non-zero function  $x \in L_M$ satisfy  the
condition $\|x\| \le 1$ and  in addition for each point  $\omega
\in \Omega$ where   $x(\omega) \ne 0$ the inequality  $|x(\omega)|
\ge h\|x\|$ hold, where $h$ is a certain positive number. Then
repeating the argument of the proof of Theorem  1 and applying
(\ref{(21)}) instead of
 (\ref{(7)}) and the monotone property of $\phi$ with respect to $u$ we get
\begin{equation}
\label{(23)}
{\cal M}(x) \ge \phi(\|x\|, h).
 \end{equation}

Now let  $x \in L_M$ be a non-zero function such that  $\|x\| \le
1$. Let
 $h$ be any number from the interval $\left(0,\D\frac{1}{\|
{\bf 1}\|}\right)$
 and
$$
x_h(\omega) := \left\{\begin{array}{ccc}
x(\omega) & \text{при} & |x(\omega)| \ge h\|x\|, \\
0 & \text{при} & |x(\omega)| < h\|x\|.
\end{array}\right.$$
 We have
$$
\|x_h\| \ge \|x\| - \|x - x_h\| \ge \|x\| - h\|  {\bf 1}\| \cdot
\|x\|=\|x\|\left(1- h\|  {\bf 1}\| \right)> 0.
$$
Since $|x(\omega )| \ge |x_h\omega )|$ it follows that (\ref{(23)}) and
the monotone property of $\phi$ with respect to $\lambda$ imply
$$
{\cal M}(x)\ge{\cal M}(x_h)\ge \phi((1 - h\|  {\bf 1}\|)\|x\|,h) . \ \ \ \ \blacksquare
$$


\begin{thebibliography}{99}

\bibitem{KR1}
    \textbf{Krasnosel'skij M.A., Rutitskij Ja.B.}
    \emph{Convex functions and Orlicz spaces.}
        Moscow: Nauka, Glavnaja redaktsija
    fiziko-matematicheskoi literatury, 1958. - 271 pp. (Russian)

\bibitem{KR2}
    \textbf{   Krasnosel'skij M.A., Rutitskij Ja.B.}
    \emph{    Orlicz spaces and     nonlinear integral equations.}
        Trudy Moskovskogo matematicheskogo obschestva,
    {\bf 7} (1958), 63-120. (Russian)

\bibitem{G}
    \textbf{   Gel'man I.V.}
    \emph{    On a certain nonlinear operator.}
        Doklady AN SSSR, {\bf 120} (1958), No 3, 454-456 (Russian).

\bibitem{R}
    \textbf{Rutitskij Ja.B.}
    \emph{On a certain property of the Orlich norm.}
     Doklady AN SSSR,
    {\bf 138} (1961), No 1, 56-58. (Russian)

\bibitem{H}
    \textbf{Honda K.}
    \emph{On a Ya.B. Ruticki's theorem concerning a property of the Orlicz norm.}
    Proc. Japan. Acad., {\bf 38} (1962), No 7, 300-303.

\bibitem{GM}
    \textbf{    Gusejnov A.I., Mukhtarov Kh.Sh.}
    \emph{    An introduction to the theory of nonlinear integral equations.}
        Moscow: Nauka, Glavnaja redaktsija
    fiziko-matematicheskoi literatury, 1980. - 416 pp. (Russian)

\bibitem{ZR}
    \textbf{   Zabreiko P.P., Rogozin S.V.}
    \emph{
    On the solvability and unique solvability of nonlinear singular integral equations
    in Oricz spaces.}
        Doklady Akademii nauk Belarusi,
    {\bf 36} (1992), No 5, 398-402. (Russian)

\end{thebibliography}
\end{document}